\documentclass[a4paper, 12pt]{article}
\usepackage{authblk}

\usepackage{import}

\usepackage[T2A]{fontenc}
\usepackage[utf8]{inputenc}
\usepackage[russian,english]{babel}

\usepackage{hyperref}
\hypersetup{
    colorlinks=true,
    urlcolor=blue,
    }

\usepackage[dvipsnames]{xcolor}

\usepackage{bm}

\usepackage{etoolbox}
\makeatletter
\def\gnewcommand{\g@star@or@long\gnew@command}
\def\grenewcommand{\g@star@or@long\grenew@command}
\def\g@star@or@long#1{% 
  \@ifstar{\let\l@ngrel@x\global#1}{\def\l@ngrel@x{\long\global}#1}}
\def\gnew@command#1{\@testopt{\@gnewcommand#1}0}
\def\@gnewcommand#1[#2]{%
  \kernel@ifnextchar [{\@gxargdef#1[#2]}%
                {\@argdef#1[#2]}}
\let\@gxargdef\@xargdef
\patchcmd{\@gxargdef}{\def}{\gdef}{}{}
\let\grenew@command\renew@command
\patchcmd{\grenew@command}{\new@command}{\gnew@command}{}{}
\makeatother

\usepackage{environ}
\usepackage{amsthm}

\usepackage{array}
\usepackage{diagbox}

\usepackage[disable]{todonotes}
\usepackage{paralist}
\newcommand{\letthengen}[4]{
        \renewcommand{\labelitemi}{ }
		#3
			\begin{inparaitem}
				#1
			\end{inparaitem}
		#4
			\begin{inparaitem}
				#2
			\end{inparaitem}
        \renewcommand{\labelitemi}{\textbullet}
}
\newcommand{\letthen}[2]{
\renewcommand{\item}{}
\letthengen{#1}{#2}{Let}{Then}
}

\makeatletter
\newenvironment{sqcases}{%
  \matrix@check\sqcases\env@sqcases
}{%
  \endarray\right.%
}
\def\env@sqcases{%
  \let\@ifnextchar\new@ifnextchar
  \left\lbrack
  \def\arraystretch{1.2}%
  \array{@{}l@{\quad}l@{}}%
}
\makeatother

\theoremstyle{plain}
\newtheorem{thm}{Theorem}
\newtheorem{lm}{Lemma}
\newtheorem{pr}{Proposition}
\newtheorem{st}{Statement}
\newtheorem{corollary}{Corollary}

\NewEnviron{subproof}[1][Subproof]
{
        
        \begin{proof}[#1]
            \BODY
        \end{proof}
}

\theoremstyle{definition}
\newtheorem{definition}{Def}

\theoremstyle{remark}
\newtheorem*{remark}{Remark}

\theoremstyle{remark}
\newtheorem{example}{Example}

\NewEnviron{EqLN}
{
        \begin{equation*}
            \begin{gathered}
                \BODY
            \end{gathered}
        \end{equation*}
}
\NewEnviron{EqLNum}
{
        \begin{equation}
            \begin{gathered}
                \BODY
            \end{gathered}
        \end{equation}
}

\makeatletter
\newcommand{\ostar}{\mathbin{\mathpalette\make@circled\star}}
\newcommand{\make@circled}[2]{%
  \ooalign{$\m@th#1\smallbigcirc{#1}$\cr\hidewidth$\m@th#1#2$\hidewidth\cr}%
}
\newcommand{\smallbigcirc}[1]{%
  \vcenter{\hbox{\scalebox{0.77778}{$\m@th#1\bigcirc$}}}%
}
\makeatother

\usepackage{relsize}

\usepackage{extarrows}

\usepackage{amssymb}

\newcommand{\ints}{\mathbb{Z}}

\newcommand{\cmplx}{\mathbb{C}}

\newcommand{\hei}{\bm{H}}

\newcommand{\supp}{\mathrm{supp} \:}
\usepackage{centernot}% http://ctan.org/pkg/centernot
\usepackage[normalem]{ulem}
\newcommand{\Sum}{\displaystyle \sum}
\newcommand{\Oplus}{\displaystyle \bigoplus}
\newcommand{\CCup}{\displaystyle \bigcup}
\newcommand{\SqCup}{\displaystyle \bigsqcup}
\newcommand{\inv}{^{-1}}

\newcommand{\Der}{\mathrm{Der}}
\newcommand{\IDer}{\mathrm{InnerDer}}
\newcommand{\ODer}{\mathrm{OuterDer}}

\newcommand{\ZDer}{\mathrm{ZDer}}
\renewcommand{\phi}{\varphi}
\newcommand{\subgrdgen}[1]{\Gamma_{#1}}
\newcommand{\subgrd}[1]{\subgrdgen{[#1]}}
\newcommand{\Obj}{Obj}
\newcommand{\Hom}{Hom}
\newcommand{\gralg}{\cmplx[G]}
\newcommand{\chcmm}[2]{ \{ #1, #2 \} }
\newcommand{\src}[1]{\mathrm{S}(#1)}
\newcommand{\trgt}[1]{\mathrm{T}(#1)}

\newcommand{\biggset}[1]{{\Bigg \{} #1 {\Bigg \}}}
\newcommand{\de}{\partial}

\usepackage{mathtools} % https://tex.stackexchange.com/a/118217/222164

\usepackage{geometry}
\usepackage{seqsplit}

\usepackage[inline]{enumitem}   
\makeatletter
\newcommand{\inlineitem}[1][]{%
	\ifnum\enit@type=\tw@
	{\descriptionlabel{#1}}
	\hspace{\labelsep}%
	\else
	\ifnum\enit@type=\z@
	\refstepcounter{\@listctr}\fi
	\quad\@itemlabel\hspace{\labelsep}%
	\fi}
\makeatother
\def\signed #1 (#2){{\unskip\nobreak\hfil\penalty50
		\hskip2em\hbox{}\nobreak\hfil\sl#1\/ \rm(#2)
		\parfillskip=0pt \finalhyphendemerits=0 \par}}

\usepackage{graphicx}
\graphicspath{../imgs/}
\usepackage{caption}

\usepackage[final]{pdfpages}

\newcommand{\comment}[1]{}

\usepackage{cleveref}

\newcommand{\Bcr}[1]{\textbf{\Cref{#1}}}

\crefname{lm}{lemma}{lemmas}
\Crefname{lm}{Lemma}{Lemmas}
\crefname{thm}{theorem}{theorems}
\Crefname{thm}{Theorem}{Theorems}
\crefname{pr}{proposition}{propositions}
\crefname{pr}{Proposition}{Propositions}
\crefname{definition}{def}{defs}
\crefname{definition}{Def}{Defs}
\crefname{corollary}{corollary}{corollaries}
\Crefname{corollary}{Сorollary}{Сorollaries}
\crefname{EqLN}{}{}
\crefname{st}{statement}{statements}
\Crefname{st}{Statement}{Statements}

\usepackage{thmtools}
\usepackage{thm-restate}

\title{Grading Structure for Derivations of Group Algebras}
\author{Andronick Arutyunov \and Igor Zhiltsov}
\date{X-XII.2022}

\begin{document}

\maketitle

\begin{abstract}
    In this paper we give a way of equipping the derivation algebra of a group algebra with the structure of a graded algebra. The derived group is used as the grading group. For the proof, the identification of the derivation with the characters of the adjoint action groupoid is used. These results also allow us to obtain the analogous structure of a graded algebra for outer derivations. A non-trivial graduation is obtained for all groups that are not perfect.
\end{abstract}

Calculation of derivations in a group algebra is a well-known problem. Present work elaborates results from articles \cite{AruMisSht16,AleAru20,Aru20} focused on studying derivations in terms of characters of adjoint action gruppoid. An important result of this research are handy formulas for quick calculation of derivations. These articles explore derivations' link to combinatorial properties of the group.

%Исследование дифференцирований в групповых алгебрах -- хорошо известная задача. Данна работа является развитием результатов ранее полученных в работах ... посвященных исследованию дифференцирований в терминах характеров группоида присоединенного действия. Одним из важных приложений является возможность получения явных формул для дифференцирований, удобных для вычислений. Данные формулы выявляют связь дифференцирований с комбинаторными свойствами группы.

Among applications note use in coding theory (see \cite{AruKos2021,CREEDON2019247})\todo{arxiv link}, Novikov algebras (see recent work \cite{kolesnikov2023prenovikov}) and more general constructions, like $(\sigma,\tau)-$derivations (see \cite{alekseev2020sigmatauderivations}).

%Среди приложений отметим приложение к некоторым задачам криптографии (см. Кридон). Приложения к данной задаче были исследованы в (https://doi.org/10.1016/j.ffa.2021.101921).  Другое актуальное направление исследований в котором могут быть полезны явные формулы для дифференцирований -- исследование алгебр Новикова (см. недавнюю работу https://arxiv.org/abs/2305.07371). Отметим, что данные результаты могут быть обобщены и на другие, более общие, алгебры. Например в https://arxiv.org/abs/2008.00390 были изучены приложения к $(\sigma,\tau)-$derivations. 

%Цель настоящей работы воспользовавшись отождествлением дифференцирований с характерами на категориях (все необходимые определения мы дадим ниже) наделить алгебру дифференцирований структурой градуированной алгебры. Основным результатом будет следующая теорема

Aim of the present work is grading the derivation algebra by identifying derivations and characters on a certain groupoid (all the necessary definitions are given in \Bcr{sect:preliminaries}). The main result of the paper follows.

Let $N$ be a fixed normal subgroup in $G$ such that $G/N$ is abelian.

\newcommand{\victim}{N}
\newcommand{\kgggen}[1]{#1_{k \in G/\victim}}
\newcommand{\kggdsum}{\kgggen{\Oplus}}
\newcommand{\kggsum}{\kgggen{\Sum}}
\newcommand{\sder}{\kggsum \Der_k}
\newcommand{\dsder}{\kggdsum \Der_k}
\newcommand{\dsderprop}{\forall k, l \in G/\victim : [\Der_k, \Der_l] \subset \Der_{kl} }

\begin{restatable}{thm}{grad}\label{grad}
    If $|G/\victim| > 1$, $\Der$ is graded with $G/\victim$, that is
    \begin{EqLN}
        \Der = \dsder, \\
        \dsderprop.
    \end{EqLN}
\end{restatable}

Here $Der_k$ is a subalgebra of derivation whose characters' support is localised entirely in one coset $a\victim = k$.

% Формулировка теоремы о градуировках. .

The structure of the work follows. \Bcr{sect:preliminaries} provides main definitions and propositions. \Bcr{sect:constr} describes the construction of grading and contains the main result and its proof. \Bcr{sect:ex} provides an example of grading for $G$ equal to discrete Heisenberg group along with an example of localising central derivations for such groups $G$ that $G$ is not a stem group. 

\todo{T=translation, R=ref, D=def, W=wording, E=elaborate}

Fix an infinite finitely-generated group $G$ for the rest of the text. 
\todo{ADD ASSUMPTIONS!}
\todo{LE TO SUBSET}
\todo{CHECK LET CASES!}
\newcommand{\FinGenAssumptionReminder}{\footnote{$G$ is finitely generated by assumptions in the beginning of the text.}}

%\section{Introduction}

\section{Preliminaries}\label{sect:preliminaries}

Recall that  \textbf{Group algebra} $\gralg$ is an algebra of formal finite sums of type ($a_1, \dots, a_n \in \cmplx, g_1,\dots,g_n \in G$)
        \begin{EqLN}
            a_1 g_1 + \dots + a_n g_n
        \end{EqLN}

%\begin{definition}
%    A \textbf{linear operator} is a function $X: \gralg \to \gralg$, such that for all $\alpha, \beta \in \cmplx \text{ и } a, b \in \gralg$
%    \begin{EqLN}
%       X(\alpha a + \beta b) = \alpha X(a) + \beta X(b)
%    \end{EqLN}
%\end{definition}

We define \textbf{derivation} as a linear operator $d$ that satisfies the Leibniz rule (for all $a, b \in \gralg$)
    \begin{EqLN}
        d(ab) = d(a) \cdot b + a\cdot d(b)
    \end{EqLN}

%\begin{definition}
%    \textbf{Commutator} in $\gralg$ is given by ($a, b \in \gralg$)
%    \begin{EqLN}
%        [a, b] := ab - ba
%    \end{EqLN}
%    \textbf{Commutator of two derivations} $d, f$ is given by ($x \in \gralg$)
%    \begin{EqLN}
%        [d, f](x) := [d(x), f(x)] = d(f x) - f(d x)
%    \end{EqLN}
%    \textit{We will refer to \textbf{commutator of two derivations} just as \textbf{commutator}.  }
%\end{definition}

Derivations over a group algebra form a Lie algebra with respect to commutator. \todo{R}

    \textit{We will denote this algebra as $\Der$ or $\Der(\cmplx[G])$.}

\subsection{Characters}

We use the technique of characters following \cite{AruMisSht16,AleAru20,Aru20}.

\newcommand{\uvcompstrip}{v_2 u_1, v_2 v_1}
\newcommand{\uvcomp}{(\uvcompstrip)}

\begin{definition}\label{md:def:grd}
    For a given group $G$ consider a \textbf{small groupoid} $\Gamma$:
    \begin{enumerate}
        \item \textbf{objects} ($\Obj$) --- elements of $G$,
        \item \textbf{arrows}\todo{R a more intuitive definition} ($\Hom$) --- pairs of elements of $G$. For an arrow $(u, v)$ its source $\src{u, v}$ is given by $v\inv u$, and its target $\trgt{u, v}$ --- by $uv\inv$ ($\Hom(a, b)$ denotes a set\todo{set?} of all arrows for which the source is $a$ and target is $b$), 
        \item\label{md:def:it:comp} Consider two arrows $\phi = (u_2, v_2) \in \Hom(b, c), \psi = (u_1, v_1) \in \Hom(a, b)$ (we will call a(n ordered) pair of arrows $\phi, \psi$ such that $\src \phi = \trgt \psi$ \textbf{composable}). The \textbf{composition} for these two arrows is given by:
        \begin{EqLN}\label{cmps}
            (u_2, v_2) \circ (u_1, v_1) := \uvcomp
        \end{EqLN}
    \end{enumerate}
\end{definition}

%\begin{st}
%    For arrows $(u_2, v_2) \in \Hom(b, c), (u_1, v_1) \in \Hom(a, b)$
%    \begin{EqLN}
%        (u_2, v_2) \circ (u_1, v_1) = \uvcomp \in \Hom(a,c)
%    \end{EqLN}
%\end{st}
%\begin{proof}
%    \begin{itemize}
%        \item Check $\src{\uvcompstrip} = a$:
%        \begin{EqLN}
%            \src{\uvcompstrip} = (v_2 v_1)\inv \cdot (v_2 u_1) = v_1\inv v_2\inv \cdot v_2 u_1 = \\ = v_1\inv u_1 =  \src{u_1, v_1} = a
%        \end{EqLN}
%        \item Check $\trgt{\uvcompstrip} = c$:
%        \begin{EqLN}
%            \trgt{\uvcompstrip} = (v_2 u_1) \cdot (v_2 v_1)\inv =  v_2 (u_1 v_1\inv) v_2\inv = v_2 \trgt{u_1,v_1} v_2\inv  = \\ = v_2 \src{u_2, v_2} v_2\inv = v_2 v_2\inv u_2 v_2\inv = u_2 v_2\inv = \trgt{u_2,v_2} c
%        \end{EqLN}
        
%    \end{itemize}
%\end{proof}

    The formula for composition does not comprise $u_2$ since if we have a pair of composable arrows $(u_2, v_2), (u_1, v_1)$, $u_2$ can be expressed in terms of $u_1, v_1, v_2$. The reader may consider this as an exercise. $\Gamma$ is the groupoid of group's inner action on itself.\todo{R}

Fix an element $a$ of $G$. Define following symbols: 
\begin{definition}
    \begin{itemize}
        \item $[a] = \{ xax\inv : x \in G \}$ is \textbf{$a$'s conjugacy class} in $G$,
        \item \begin{EqLN}
            G^G := \{[g] : g \in G \}
        \end{EqLN}
        \item $\subgrd a$ is $\Gamma$'s subgroupoid, \textit{informally, a connected component in $\Gamma$}, given by:
        \begin{EqLN}
            \Obj(\subgrd a) := [a] = \{x \in Obj : x \in [a] \} \\
            \Hom(\subgrd a) := \{ (u, v) \in \Hom : u,v\in[a]\}
        \end{EqLN}
    \end{itemize}
\end{definition}

\begin{definition}
    A \textbf{character} on $\Gamma$ is a function $\chi: \Hom \to \cmplx$, such that:
    \begin{itemize}
        \item \textbf{\textit{(Composition)}} for each pair of composable arrows $\phi, \psi$:
            \begin{EqLN}
                \chi(\phi \circ \psi) = \chi(\phi) + \chi(\psi)
            \end{EqLN}
        \item \textbf{\textit{(Locally finite)}} $\forall \bm y \in G$ there is a \textbf{finite} set of $\bm x \in G$, such that $\chi(x, y) \neq 0$.\todo{T}
    \end{itemize}
\end{definition}

Holds the following decomposition:

\begin{lm}(Decomposition)\label{decompx}
    $\Gamma = \SqCup_{[a] \in G^G} \subgrd a$
\end{lm}\todo{D cup?}\todo{R}

\begin{remark}
    Although characters being locally finite may seem as an alien and a bit too technical detail, it is deliberately placed in the definition to stress that we will \textbf{not} consider \textit{"non-locally finite characters"}. The reasons will become clear, among the rest, in \Bcr{difffla}. 
\end{remark}

%\begin{st}\label{maindef:pr:1}
%    $\Hom(a, b)$ is non-empty $\iff$ $\exists x \in G: a = xbx\inv$. \todo{R}
%\end{st}

We will need the following statement:

\begin{st}\label{maindef:lm:1}\label{maindef:pr:1}
    \letthengen{
        \item $(u, v) = \phi \in \Hom$,
        \item $a \in G$.
    }{
        \item $\phi \in \Hom(\subgrd a)$,
        \item $\src\phi \in [a]$,
        \item $\trgt\phi \in [a]$.
    }{Let}{Then the following statements are equivalent}
\end{st}

It is proved by direct calculation.

%\begin{proof}
%    By definition, $\src\phi = v\inv u, \trgt\phi = uv\inv$. Thus:
%    \begin{itemize}
%        \item If $\phi \in \Hom(\subgrd a)$, then $\src\phi \in [a]$ by definition.\todo{W}
%        \item If $\src\phi \in [a]$ then $\trgt\phi \in [a]$ by \Bcr{maindef:pr:1}.
%        \item If $\trgt\phi \in [a]$, then by \Bcr{maindef:pr:1} $\src\phi \in [a]$, thus $\phi \in \Hom(\subgrd a)$ by definition\todo{W}.
%    \end{itemize}
%\end{proof}

\subsection{Connection between Characters and Derivations}

The following theorem motivates to consider \textbf{(locally finite)} characters when studying derivations. \textit{Informally, \Bcr{difffla} shows that characters may be seen as a generalization of linear operator's matrix.}

\begin{thm}[Derivation formula and derivation character, {\cite[section~2]{AruMisSht16}}]\label{difffla}\todo{T}
    For each derivation $d$ there exists a unique character $\chi$ such that for each $x \in G$ holds
        \begin{equation}
            d(x) = \Sum_{k \in G} \chi(k, x) k
        \end{equation}
\end{thm}

Consider $d$, $\chi$ from \Bcr{difffla}. We will say that character $\chi$ \textbf{gives} derivation $d$ (derivation $d$ \textbf{is given by character} $\chi$; we will omit words "derivation" and "character"). For a derivation $d$ let $\chi^d$ be a character such that $\chi^d$ gives $d$.  \todo{T}

\Bcr{difffla} implies:

\begin{corollary}\label{mp:corollary:1}
    \letthen{
        \item $d, \de$ be derivations given by characters $\chi^d, \chi^\de$ correspondingly,
    }{
        \item $d + \de$ be given by $\chi^d+\chi^\de$.
    }
\end{corollary}
%\begin{proof}
%    \Bcr{difffla} implies
%    \begin{EqLN}
%        d(x) = \Sum_{k \in G} \alpha(k, x) k, \quad \quad \de(x) = \Sum_{k \in G} \beta(k, x) k
%    \end{EqLN}
%    Thus
%    \begin{EqLN}
%        (d+\de)(x) = \Sum_{k \in G} (\alpha(k, x) + \beta(k, x)) k
 %   \end{EqLN}
%\end{proof}

\begin{definition}
    \letthen{
        \item $d$ be given by $\alpha$,
        \item $\de$ be given by $\beta$.
    }{
        \item $\chcmm{\alpha}{\beta}$ is the character that gives $[d, \de]$. 
    }
\end{definition}

\begin{st}[\textit{"Matrix" product}, {\cite[Proposition~2.4]{AleAru20}}]\label{mtrx}
    \letthen{
        \item $\alpha, \beta$ be characters. 
    }{
        \item $ \chcmm{\alpha}{\beta} $ satisfies ($a, b \in G$)
        \begin{EqLN}
            \chcmm{\alpha}{\beta}(a, b) = \Sum_{k \in G} \alpha(a, k) \beta(k, b) - \beta(a, k) \alpha(k, b) 
        \end{EqLN}
    }
\end{st}

Two examples of derivations follow. \Bcr{innerex} will be needed to prove \Bcr{grad}.

\newcommand{\innerdef}{Let $a \in G$. Recall that derivation $d_a$ is called \textbf{inner} if for any $x \in \gralg$
$$
    d_a(x) = [x, a] = xa - ax
$$
}

\innerdef

%Let $a \in G$. Recall that derivation $d_a$ is called \textbf{inner} if for any $x \in \gralg$
%     \begin{EqLN}
%        d_a(x) = [x, a] = xa - ax
%    \end{EqLN}

\begin{example}[Character of inner derivation {\cite[Proposition~3]{Aru20}}]\label{innerex}
    \letthen{
        \item $a\in G$.
    }
    {
        \item character $\chi_a$ given by formula
        \begin{equation}\label{innerchar}
            \chi_a(\phi) = 
                \begin{dcases}
                    \: \: \: 1,  & a = \src{\phi}, \\
                    -1, & a = \trgt{\phi}, \\
                    \: \: \: 0,  & \text{otherwise.}
                \end{dcases}   
        \end{equation}
        gives $d_a(x) = [x, a]$.
    }
\end{example}

\begin{example}\todo{central}
    Another possible example of derivations are \textit{central derivations}. We will call derivation $d$ central if there exists such central element $z \in Z(G)$ and homomorphism $\tau: G \to (\cmplx, +)$ such that for all \textit{basis elements} $g \in G$:
    \begin{EqLN}
        d(g) = \tau(g)gz
    \end{EqLN}

    Such an operator is indeed a derivation, see \cite[Proposition~4]{Aru20}. \cite[Proposition~5]{Aru20} shows that non-trivial central derivations are not inner. Moreover, \cite[Proposition~6]{Aru20} shows that central derivations form a Lie subalgebra in $\Der(G)$.
\end{example}

\begin{definition}
    %The set of arrows in $\Gamma$ at which character $\chi$ does not vanish is \textbf{$\bm \chi$'s support ($\bm{\supp \chi}$)}.
    For a given character $\chi$ we define support of $\chi$ as following
    \begin{EqLN}
        \supp \chi = \{\phi \in \Hom : \chi(\phi) \neq 0 \}
    \end{EqLN}
\end{definition}

For the given subset $M \subset G$ denote by $\bm {\Der_M}$ the set of derivations $d$ such that for character $\chi$ that gives $d$: $\supp \chi \subset M$.

\begin{example}\label{innerextwo}
    Recall character $\chi_a$ from \Cref{innerex} (where $a \in G$.) Its support is easily calculated
    \begin{EqLN}
        \supp \chi_a = \{ \phi : \src{\phi} = a \} \cup \{ \psi : \trgt{\psi} = a \} 
    \end{EqLN}

    By \Bcr{maindef:pr:1}, we can localize $\supp \chi_a$ in a single conjugacy class $a$ (we will need such a localisation later in \Bcr{grad})
    \begin{EqLN}
        \supp \chi_a \subset \subgrd a
    \end{EqLN}

\end{example}
\subsection{Applying Decomposition}\label{sub:appdecomp}

\Bcr{decompx} establishes decomposition of groupoid $\Gamma$. The current section presents \textit{decompositions for \textbf{(locally finite)} characters and derivations.}\todo{W???}

%{ \color{red}
%Это утверждение справедливо не для всех характеров, а только для локално финитных. Собственно где-то там выше надо было уточнить, что речь мы ведем именно о них. И вообще, дифференцированиям на групповой алгебре соответствуют как раз локально финитные характеры (иначе мы не попадем в групповую алгебру).
%}

\newcommand{\compgoal}{\supp \chi \le \CCup_{k=1}^N \subgrd {a_k}}
\newcommand{\myfiniteset}{\{ [u] \in G^G: \exists x \in G: d_u(x) \neq 0 \}}

The following two lemmas are equivalent. We prove the first one.

\begin{lm}\label{localfin}\todo{Name?}
    \letthen{
        \item $\chi$ be a character.
    }{
        \item there exists finitely many $a_1,\dots,a_N \in G$ such that
        \begin{EqLN}
            \compgoal
        \end{EqLN}
    }
\end{lm}
\begin{lm}[Derivation decomposition] \label{derdecomp}\label{finderdecomp}
    Then for each $ u \in G$ such that $\chi_u$ is a character, holds decomposition $d = \Sum_{[u] \in G^G} d_u$, and the set $\myfiniteset$ is finite.
\end{lm}

\begin{proof}[Proof for \Bcr{localfin}]
    Let $G = \langle X \mid R \rangle$, where $X =: \{x_1, \dots, x_k\}$ is finite \todo{D} ($G$ is assumed to be finitely-generated throughout the text). 
    
    Consider $(u, v) \in \Gamma$. Let $n = n(v)$ be minimal nonnegative integer such that
    \begin{EqLN}
        \exists y_0, \dots, y_n \in X \cup X\inv: v = y_0y_1 \dots y_n
    \end{EqLN}
    
    \begin{enumerate}
        \item\label{appdecomp:item:1} Show that
        \begin{EqLNum}\label{appdecomp:eq:1}
            \exists z_0, \dots, z_n \in G: (z_0, y_0) \circ \dots \circ (z_n, y_n) = (u,v)
        \end{EqLNum}
        \begin{subproof}
            Induction by $n = n(v)$.
            \textbf{Base:} $n=0$ --- $z_0 = u$.
            \textbf{Step:} Consider $z_0 = uv\inv y_0$. Then $(z,y_0), (y_0\inv u, y_0\inv v)$ are composable since
        \begin{EqLN}
            \src{z_0,y_0} = y_0\inv z_0 = y_0\inv uv\inv y_0 = \\ = y_0\inv u(y_0\inv v)\inv = \trgt{y_0\inv u, y_0\inv v}
        \end{EqLN}

        Moreover,
        \begin{EqLN}
            (z_0,y_0) \circ (y_0\inv u, y_0\inv v) = (u, v)
        \end{EqLN}

        Notice that $y_0\inv v = y_1\dots y_n$, thus $n(y_0\inv v) < n = n(v)$. Therefore, applying the step of induction, get \cref{appdecomp:eq:1}.
        \end{subproof}

        \item Since $\chi$ is \textbf{\textit{locally finite}}, the set $B = (G \times (X \cup X\inv)) \cap \supp \chi$ is finite. By \Bcr{decomp} for each arrow $\phi$ there exists a unique element $a \in G$ such that $\phi \in \Hom(\subgrd a)$; thus, there exists a finite set $A = \{a_1, \dots, a_N\}$ such that for any $a \notin A$: $B \cap \Hom(\subgrd a) = \O$. Thus, by \cref{appdecomp:item:1} for any $a \notin A$: $\supp \chi \cap \Hom(\subgrd a) = \O$. Therefore:
        \begin{EqLN}
            \compgoal
        \end{EqLN}
    \end{enumerate}
\end{proof}

A very nice alternative proof for \Bcr{localfin} was submitted in an anonymous review.

\begin{proof}[Alternative proof for \Bcr{localfin}]
    Let $d$ be the derivation given by character $\chi$. 
    
    %Let $S$ be a finite generating set for $G$. For each $s \in S$ consider an element $s\inv d(s)$. There exists a finite subset $M \subset G$ such that for any $s \in S$ there exist complex numbers $a_m, m \in M$ such that
    %\begin{EqLN}
    %    s\inv d(s) = \Sum_{m \in M} a_m m
   % \end{EqLN}
   % Consider a finite union of conjugacy classes 
    %\begin{EqLN}
    %    P = \CCup_{m \in M} [m] 
    %\end{EqLN}
   % \textit{The definitions will become clear later. Our goal is to show that $\supp \chi \subset \CCup_{m \in M} \subgrd m$.}

    Let $P$ be a union of conjugacy classes such that
    \begin{EqLN}
        \supp \chi \subset \CCup_{g \in P} \subgrd g =: U
    \end{EqLN}

    The following statements are equivalent:
    \begin{itemize}
        \item $(x, y) \in \Hom(U)$,
        \item $y\inv x = \src{x,y} \in P$
    \end{itemize}

    Consider an element $y \in G$ such that
    \begin{EqLN}
        y\inv d(y) \in \langle P \rangle
    \end{EqLN}

    \textit{Here $\langle X \rangle$ for $X \subset G$ denotes a set of all finite sums $a_1x_1 + \dots + a_nx_n$ such that $a_1, \dots, a_n \in \cmplx$ and  $x_1, \dots, x_n \in X$.}

    Let's calculate $y\inv d(y)$ by \Bcr{mtrx}.
    \begin{EqLN}
        y\inv d(y) = y\inv \Sum_{x \in G} \chi(x, y)x = \Sum_{x\in G} \chi(x,y) y\inv x = \Sum_{x\in G} \chi(x,y) \src{x,y} \in \langle P \rangle
    \end{EqLN}

    Therefore, for each $x$ such that $\chi(x, y) \neq 0: y\inv x \in P$.

    Consider a set
    \begin{EqLN}
        H = \Big\{ y : y\inv d(y) \in \langle P \rangle \Big \}
    \end{EqLN}

    As a simple exercise, check $H$'s being a subgroup in $G$.

    To summarize, if $y$ is in subgroup $H \le G$ (that is $y\inv d(y) \in \langle P \rangle$) then for each $x$ such that $\chi(x, y) \neq 0: y\inv x \in P$.

    To finish the proof let's choose such $P$ that $P$ is a union of a \textbf{finite number} of conjugacy classes and $H = G$. 
    
    To achieve this, consider a finite generating set $S$ for $G$ and a \textbf{finite} subset $M \subset G$ such that for any $s \in S$ there exist complex numbers $a_m, m \in M$ such that
    \begin{EqLN}
        s\inv d(s) = \Sum_{m \in M} a_m m
   \end{EqLN}
   \textit{Informally, calculate all $s\inv d(s), s \in S$, which are finite sums, and store all elements in $G$ present in at least one of finite sums.}

   Since $M$ is finite, $P = \CCup_{m \in M} [m]$ is a union of finite number of conjugacy classes. Moreover, for such $P$
   \begin{EqLN}
       S \subset H
   \end{EqLN}
   Thus, $G = H$. 

   All in all, there exists a union $P$ of a \textbf{finite} number of conjugacy classes such that
   \begin{EqLN}
        \supp \chi \subset \CCup_{g \in P} \subgrd g =: U
    \end{EqLN}
\end{proof}

 Let $d$ be the derivation given by character $\chi$, and
         $$\chi_u(\phi) := 
        \begin{cases}
            \chi(\phi), &\phi \in \Hom(\subgrd u), \\
            0,          &otherwise
        \end{cases}.$$
We will denote the derivation given by $\chi_u$ as $d_u$.

\newcommand{\dudescr}{$d_u$ --- дифференцирование, характер которого совпадает с характером $d$ на $\subgrd u$ и равен 0 вне него; почти все $d_u$ нулевые}

\section{Constructing Graded Algebra}\label{sect:constr}

%\import{sections/}{CommChar.tex}

\subsection*{Grading with Abelian Quotients}\label{sect:grading}

\begin{definition}\label{grad:def:1}
    \letthen
    {
    \item $A$ be an abelian group,
    \item $e$ be a neutral element in $A$,
    \item $\mathfrak{A}$ be a Lie algebra,
    \item $\mathfrak A$ can be expressed as a direct sum
    \begin{EqLN}\label{grading:def:grading}
        \mathfrak{A} = \Oplus_{n \in A} \mathfrak{A}_n, \text{ such that } \\
        \forall n, l \in A: [\mathfrak{A}_n, \mathfrak{A}_l] \subset \mathfrak{A}_{nl}
    \end{EqLN}
    such that $\mathfrak{A}_e \neq \mathfrak{A}$.
    }
    {
    \item $\mathfrak A$ is called \textbf{graded (with $\bm A$)}. The direct sum \cref{grading:def:grading} is called \textbf{$\bm{\mathfrak{A}}$'s grading with $\bm A$}.
    }
\end{definition}

Notice that trivial gradings are \textbf{excluded}.

%The main goal of \Bcr{sect:grading} is to prove that $\Der$ is graded.

\begin{definition}
    A commutator subgroup of group $G$ is
    \begin{EqLN}
        G' := \{ x y x\inv y\inv : x,y \in G \}
    \end{EqLN}
\end{definition}

Recall a few well-known definitions and statements that will be required below.

\begin{st}\label{commnorm}
    For any group $G$ a subgroup $G'$ is normal,
        $G' \triangleleft G$
\end{st}

Let $N$ be a fixed normal subgroup in $G$ such that $G/N$ is abelian. \textit{Conceptually, the latter condition derives from the following \Bcr{grad:def:1}, however it is also needed for more technical, yet crucial details, like \Bcr{normadj}.}

\begin{st}\label{commminprop}
    %For any $N \triangleleft G$: $G/N$ is abelian \textbf{iff}
    $G' \subset N$.
\end{st}

\begin{lm}\label{normadj}\label{abeladj}
\letthen{
    %\item $N \triangleleft G$ --- $N$ is a normal subgroup in $G$,\todo{D normal, triangle!}
    %\item $G/N$ is abelian,
    \item $a \in N$.
}{
    \item $[a] \subset aN$. \todo{"Это что?"}
}
\end{lm}

\begin{proof}
    A calculation for any $a, t \in G$ proves lemma:
    \begin{EqLN}
        tat\inv N = tat\inv a\inv N a = Na = aN
    \end{EqLN}
    The first and the last equivalences hold since $N$ is normal; the second equivalence holds by \Bcr{commminprop}
    since $G/N$ is abelian ($tat\inv a\inv \in G' \subset N$.)
\end{proof}  

\newcommand{\defsubgrdcl}{
    \subgrdgen{aN} = \CCup_{k \in aN} \subgrd{k}, % N \triangleleft G,% G/N \text{ is abelian} 
    \\
    \Der_{aN} = \biggset{d \in Der : \supp \chi^d \subset \subgrdgen{aN}}.
}

\textit{Note that \Bcr{normadj} would have been false if we did not require $G/N$ to be abelian. Consider $G = S_4, N=V_4, a=(12)$ for counterexample.}\todo{Delete?}

\Bcr{normadj} motivates the following symbols:
\begin{EqLN}
    \defsubgrdcl
\end{EqLN}

%\begin{corollary}[From \Bcr{normadj,commnorm}]\label{abeladj}
%    $\forall a \in G': [a] \subset aG'$. 
%\end{corollary}

\newcommand{\commsupp}[2]{\subgrd #1 \otimes \subgrd #2}
\newcommand{\commsuppab}{\commsupp{a}{b}}
\newcommand{\commsuppRHS}[4]{\CCup_{#1 \in [#3]} \CCup_{#2 \in [#4]} \subgrd{#1 #2}}
\newcommand{\commsuppRHSuvab}{\commsuppRHS{u}{v}{a}{b}}

\newcommand{\cla}{a\victim}
\newcommand{\clb}{b\victim}
\newcommand{\clab}{ab\victim}
\newcommand{\sga}{\subgrdgen\cla}
\newcommand{\sgb}{\subgrdgen\clb}
\newcommand{\sgab}{\subgrdgen\clab}
\newcommand{\homa}{\Hom(\sga)}
\newcommand{\homb}{\Hom(\sgb)}
\newcommand{\homab}{\Hom(\sgab)}

\newcommand{\commchartemplateR}[5]{\Sum_{#5 \in G} #1(#3, #5) #2(#5, #4) - #2(#3, #5) #1(#5, #4)}
\newcommand{\commchartemplate}[5]{\{#1, #2\}(#3, #4) = \commchartemplateR{#1}{#2}{#3}{#4}{#5}}
\newcommand{\commchar}{\commchartemplate{\alpha}{\beta}{h}{g}{k}}

\begin{lm}\label{commcomm}
    \letthen{
        \item $a, b \in G$, 
        \item $d$ is given by character $\alpha: \supp \alpha \subset \sga$, 
        \item $\de$ is given by character $\beta: \supp \beta \subset \sgb$.
    }{
        \item $\supp \chcmm{\alpha}{\beta} \subset \sgab$.
    }
\end{lm}

\begin{proof}

    \Bcr{mtrx} implies:
    \begin{EqLN}
        \commchar
    \end{EqLN}
    Consider an arrow $(h, g): \chcmm{\alpha}{\beta}(h, g) \neq 0$. There exists $k \in G$:
    \begin{EqLN}
        \begin{sqcases}
            \alpha(h, k)\beta(k, g) \neq 0 \\
            \beta(h, k)\alpha(k, g) \neq 0
        \end{sqcases}
    \end{EqLN}\todo{D sqcases?}

    \begin{EqLNum}\label{commchar:eq:1}
        \begin{sqcases}
            \begin{dcases}
                \alpha(h, k) \neq 0 \\
                \beta(k, g) \neq 0
            \end{dcases}
              \\ \\
             \begin{dcases}
                 \beta(h, k) \neq 0 \\
                 \alpha(k, g) \neq 0
             \end{dcases}
            
        \end{sqcases}
    \end{EqLNum}
    
    Expressing \cref{commchar:eq:1} in terms of $\supp$:
    
    \begin{EqLNum}\label{commchar:eq:2}
    \begin{sqcases}
            \begin{dcases}
                (h, k) \in \supp \alpha \subset \homa \\
                (k, g) \in \supp \beta \subset \homb
            \end{dcases}
            \\ \\
            \begin{dcases}
                (k, g) \in \supp \alpha \subset \homa \\
                (h, k) \in \supp \beta \subset \homb
            \end{dcases}
            
    \end{sqcases}
    \end{EqLNum}

    By \Bcr{maindef:lm:1}, an arrow $\phi$ belongs to $\subgrd x \subset \subgrdgen{x\victim}$ iff its target $\trgt \phi$ belongs to $[x] \subset x\victim$. Thus, \cref{commchar:eq:2} implies:
    
    \begin{EqLNum}\label{commchar:eq:3}
        \begin{sqcases}
            \begin{dcases}
                h k\inv = u \in \cla \\
                k g\inv = v \in \clb
            \end{dcases}
            \\ \\
            \begin{dcases}
                k g\inv = u \in \cla \\
                h k\inv = v \in \clb
            \end{dcases}
        \end{sqcases}
    \end{EqLNum}

    Multiplying the equations in \cref{commchar:eq:3}:

    \newcommand{\uvdef}{,\textbf{ for } u \in \cla, v \in \clb}
    
    \begin{EqLN}
        \begin{sqcases}
           h g\inv = uv \\
           h g\inv = vu
        \end{sqcases}
        \uvdef
    \end{EqLN}

    By definition,\todo{W} $\trgt{h, g} = hg\inv$. By \Bcr{maindef:lm:1}:

    \begin{EqLN}
        \begin{sqcases}
            (h, g) \in \Hom(\subgrd {uv}) \\
            (h, g) \in \Hom(\subgrd {vu})
        \end{sqcases}
    \end{EqLN}

    Since $[uv] = [u \cdot vu \cdot u\inv] = [vu]$ and $uv \in ab\victim$, then by \Bcr{abeladj}:

    \begin{EqLN}
         [uv] \subset uv\victim = ab\victim
    \end{EqLN}

    Thus:

    \begin{EqLN}
        (h, g) \in \subgrd {uv} \subset \sgab \uvdef
    \end{EqLN}

    All in all, $\chcmm{\alpha}{\beta}(h, g) \neq 0 \Rightarrow (h,g) \in \sgab$, therefore $\supp \chcmm{\alpha}{\beta} \subset \sgab$.
\end{proof}

%{\color{red}
%неудачные обозначения. $A$ нужно заменить на $Der_k$ например, и дать определение $Der_k$ до формулировки теоремы. Сама формулировка теоремы должна быть в духе: алгебра дифференцирований -- алгебра градуированная при помощи фактор группы по коммутанту. И в доказательстве первый пункт -- вообще-то самый главный, его надо прояснить.
%}

%\newcommand{\kgggen}[1]{#1_{k \in G/G'}}
%\newcommand{\kggdsum}{\kgggen{\Oplus}}
%\newcommand{\kggsum}{\kgggen{\Sum}}
%\newcommand{\sder}{\kggsum \Der_k}
%\newcommand{\dsder}{\kggdsum \Der_k}
%\newcommand{\dsderprop}{\forall k, l \in G/G' : [\Der_k, \Der_l] \subset \Der_{kl} }

%Denote set of such derivations $d$ given by such character $\chi$ that $\supp \chi \le \subgrdgen{aG'}$ as $\Der_{aG'}$.

\grad*

\begin{proof} 
Consider the sum $\sder$.

\begin{enumerate}
    \item First, show that $\sder$\todo{"Сумма прямая"; да, но ещё не доказано} is equal to $\Der$.
    \begin{itemize}
        \item $\sder \subset \Der$ --- since $\Der$ is closed under finite sums.
        \item Now we show an opposite inclusion\todo{T}. 
        Consider an arbitrary $d \in \Der$. \Bcr{derdecomp,abeladj} imply (see \Bcr{derdecomp} for definition of $d_u$)
        \begin{EqLN}
            d = \Sum_{[u] \subset G} d_u = \kggsum \Big( \Sum_{[u] \subset k} d_u \Big)
        \end{EqLN}

        Consider for given $k \in G/\victim$
        \begin{EqLN}
            s_k := \Sum_{[u] \subset k} d_u \in \Der_k
        \end{EqLN}

        By \Bcr{derdecomp}, there is only a finite number of such $[u] \subset G$ that $d_u$ is not constant 0\todo{W}. Thus, there exists an integer $N$ and $k_1, \dots, k_N \in G/\victim$ such that for any $k \in G/\victim, k \neq k_1, \dots, k_N$: $s_k$ is constant 0.

        Thus,
        \begin{EqLN}
            d = \Sum^N_{i=1} s_{k_i}
        \end{EqLN}

        Thus,
        \begin{EqLN}
            \Der \subset \sder
        \end{EqLN}
        \item All in all, $\Der = \sder$.
    \end{itemize}

    \item $\sder$ is direct.
        \begin{subproof}
            Let $d_k \in \Der_k$ be given by $\chi_k$. Suppose that
            \begin{EqLN}
                \kggsum d_k \equiv 0
            \end{EqLN}\todo{D equiv or change?}

            Since constant 0 is a derivation given by constant character (equal to 0)\todo{W}, by \Bcr{mp:corollary:1} for any $\phi \in \Hom$
            \begin{EqLN}
                \kggsum \chi_k(\phi) = 0
            \end{EqLN}

            Since $\supp \chi_k, \supp \chi_l$ are disjoint (for $k \neq l$), for each $\phi \in \Hom$ exists no more than one $k \in G/\victim$ such that $\chi_k(\phi) \neq 0$. Thus, for each $\phi \in \Hom$ and for each $k \in G/\victim: \chi_k(\phi) = 0$, and for each $k \in G/\victim: d_k \equiv 0$.\todo{D equiv? Change?}  Therefore, $\sder = \dsder$ is direct by definition of direct sum.\todo{Check}
        \end{subproof}
        
        \item We established that
        \begin{EqLN}
            \Der = \dsder
        \end{EqLN}

        Since there exists such $m \in G$ that $m\victim \neq \victim$. Thus, by \Cref{innerex,innerextwo} the inner derivation $[x, m]$ is given by a character $\chi_m$ such that $\supp \chi_m \subset [m] \subset m\victim$ by \Bcr{abeladj}. Thus, $\Der \neq \Der_{\victim}$ \textit{(i.e. the grading is not trivial.)}
        
        Finally, check that ($\forall k, l \in G/\victim$)
        \begin{EqLN}\label{maingradproperty}
            \dsderprop
        \end{EqLN}
        
        \begin{subproof}[Proof for \cref{maingradproperty}]
            Let $d_k \in \Der_k, d_l \in \Der_l$. Let character $\chi_k$ give $d_k$, character $\chi_l$ give $d_l$. \todo{W}Thus, by definition of $\Der_k, \Der_l$
            \begin{EqLN}
                \supp \chi_k \le \subgrdgen k, \quad\quad \supp \chi_l \le \subgrdgen l            \end{EqLN}

            Therefore, by \Bcr{commcomm}
            \begin{EqLN}
                \supp \chcmm{\chi_k}{\chi_l} \le \subgrdgen {kl} 
            \end{EqLN}
            Finally,
            \begin{EqLN}
                [d_k, d_l] \in \Der_{kl}
            \end{EqLN}
        \end{subproof}
    \end{enumerate}
\end{proof}
\todo{RETURN SUBSECTION}
%\import{./}{GradingOtherQ.tex}

\begin{example}
    Let $G$ be a perfect group ($|G/G'| = 1$), that is $G' = G$. \Bcr{grad} yields \textit{a trivial grading}(which we do not regard as a grading in this text) for $\Der(\cmplx[G])$ since $|G/G'| = 1$. 
\end{example}

And vice versa:

\begin{corollary}
    If $G$ is NOT a perfect group ($G \neq G'$) then $\Der$ admits a (non-trivial) grading with $G/G'$.
\end{corollary}

\begin{example}
    Let $G$ be a knot group (i.e. let there exist some knot $K$ such that $G$ is the knot group of $K$). It is well-known that in this case $G/G' = \ints$, therefore $Der(G)$ admits a grading with $\ints$.
\end{example}

%\subsection*{Grading with Other Abelian Quotients}
%\import{sections/}{GradingOtherQ.tex}

\section{Examples}\label{sect:ex}

%\subsection{Superalgebras}

%\import{Examples}{Super.tex}

\subsection*{Discrete Heisenberg Group}\todo{change to Out}

Consider discrete Heisenberg group (a group of $3 \times 3$ upper unitriangular matrices with integer entries). Following \cite{Aru20}, we use this group as a handy example since it admits easy calculations. 

\newcommand{\dhe}{\Der(\hei)}
\newcommand{\he}[3]{
\begin{pmatrix}
        1 & #1 & #3 \\
        0 & 1 & #2 \\
        0 & 0 & 1
\end{pmatrix}
}
\newcommand{\hab}{\hei / \hei'}
\newcommand{\intsq}{\ints \oplus \ints}

\begin{definition}
    Consider a group of integer unitriangular matrices with respect to matrix multiplication:
    \begin{EqLN}
        \hei = \biggset{\he{a}{b}{c}  : a, b, c \in \ints}
    \end{EqLN}
    \begin{EqLN}\label{hmul}
        \he{a}{b}{c} \he{x}{y}{z} := \he{a+x}{b+y}{c + z + ay}
    \end{EqLN}

    Since all the matrices in $\hei$ have determinant 1, the inverse is well-defined and given by
    \begin{EqLN}
        \he{a}{b}{c}\inv = \he{-a}{-b}{ab - c}
    \end{EqLN}
\end{definition}

Our goal is to grade \todo{T} $\dhe$.

\begin{definition}
    The centre of $G$ is
    \begin{EqLN}
        Z(G) = \{ z \in G: \forall g \in G: gz = zg \}
    \end{EqLN}
\end{definition}

The follwing statements are well-known and trivial.

\begin{st}\todo{R}\label{heicommsub}
    \begin{EqLN}
        \hei' = Z(\hei) = \biggset{\he{0}{0}{a} : a \in \ints}
    \end{EqLN}
\end{st}

\begin{st}\label{heiab}
    \begin{EqLN}
        \hab \simeq \intsq
    \end{EqLN}
\end{st}
%\begin{proof}
%    Consider a matrix
%    \newcommand{\ma}{\he{a}{b}{c}}
%    \begin{EqLN}
%        A = \ma \\
%    \end{EqLN}
%
%    Consider a set $A\hei'$. We will show that
%    \begin{EqLN}
%        A\hei' 
%        \xlongequal{\Bcr{heicommsub}}
%        \biggset{A\he{0}{0}{n} : n \in \ints}
%        = \\ \xlongequal{\cref{hmul}}
%        \biggset{\he{a}{b}{c + n} : n \in \ints} 
 %       = 
 %       \biggset{\he{a}{b}{c + n} : n \in \ints} 
 %       = \\ =
 %       \biggset{\he{a}{b}{c'} : c' \in \ints} 
 %   \end{EqLN}

 %   Thus, there is a well-defined bijection $\phi$ from the quotient group $\hab$ to $\intsq$ given by
 %   \begin{EqLN}
 %       \he{a}{b}{c}\hei' \xmapsto{\quad \mathllarger \phi \quad } (a, b)
 %   \end{EqLN}

%    \newcommand{\apphi}[1]{\mathllarger{\phi}{\Bigg (} #1 {\Bigg )}}
%    Finally, by direct computation we can see $\phi$ being a homomorphism %since
%    \begin{EqLN}
%        \apphi{\he{a}{b}{c} \hei' \cdot \he{x}{y}{z} \hei'} 
 %        =
    %    \apphi{\biggset{\he{a}{b}{c'} : c' \in \ints}
    %    \cdot 
    %    \biggset{\he{x}{y}{z'} : z' \in \ints}}
    %    = \\ =
    %    \apphi {\biggset{\he{a}{b}{c'}\he{x}{y}{z'} : c', z' \in \ints} }
    %    = 
    %    \apphi {\biggset{\he{a+x}{b+y}{z' + c' + ay} : c', z' \in \ints}}
    %    = \\ =
    %    \apphi {\biggset{\he{a+x}{b+y}{n} : n \in \ints}} = \apphi{\he{a+x}{b+y}{0}\hei'} 
    %    = \\ = 
    %    (a+x, b+y) = (a,b) + (x, y)
    %    =
%        \apphi{\he{a}{b}{c} \hei'} + \apphi{\he{x}{y}{z} \hei'}
%    \end{EqLN}
%
%    All in all, $\phi : \hab \to \intsq$ is an isomorphism. Therefore
%    \begin{EqLN}
%        \hab \simeq \intsq
%    \end{EqLN}
%\end{proof}

\newcommand{\psidef}{$\psi: \intsq \to \hei / \hei'$ be an isomorphism}

Let \psidef. Recall the symbols: 
\begin{EqLN}
    \defsubgrdcl    
\end{EqLN}

Define:

\newcommand{\dij}[1]{\Der_{(#1)}}

\begin{EqLN}
    \dij{i,j} := \Der_{\psi(i,j)}
\end{EqLN}

\begin{corollary}[From \Bcr{heiab}, \Bcr{grad}]\label{heigrad}
    $\dhe$\todo{D} is graded with $\intsq$, that is
    \begin{EqLN}
        \dhe = \Oplus_{(i, j) \in \intsq} \dij{i,j} \\
        \forall (i,j), (k,l) \in \intsq: [\dij{i,j}, \dij{k,l}] \subset \dij{i+k,j+l}
    \end{EqLN}    
\end{corollary}

\begin{example}\label{hei:cor:3} \Bcr{heicommsub} implies:
    if $d$ is given by such $\chi$ that $\supp \chi \le \subgrd z$ for $z \in Z(\hei)$, then $d \in \Der_{(0,0)}$. 
\end{example}

%To elaborate this structure, we take a closer look at the conjugacy classes of $\hei$.

%\begin{lm}\todo{R or prove}
%    \letthen{
%        \item $A = \he{a}{b}{c}$
%        \item $B = \he{x}{y}{z}$
%    }
%    {
%        \item 
%        \begin{EqLN}
%            [A] = [B] \iff 
%            \begin{dcases}
%               a = x \\
%                b = y \\
%                c - z \equiv 0 \mod \gcd(a,b)
%            \end{dcases}
%        \end{EqLN}\todo{D equiv, gcd}
%    }
%\end{lm}
%\begin{lm}
%    \letthen{
%        \item $(i,j) \in \intsq$
%        \item \psidef,
%        \item $a\hei' = \psi(i,j)$,
%    }
%    {
%        \item if $i^2 + j^2 > 0$:
%        \begin{itemize}
%            \item Let $N = \gcd(i,j)$,
            %\item there exist $x_1, \dots, x_N \in a\hei'$ %such that:
%        \begin{EqLN}
%            a\hei' = [x_1] \sqcup \dots \sqcup [x_N] 
%        \end{EqLN} 
%        \end{itemize}
%        \item else:
%        \begin{itemize}
%            \item
%            \begin{EqLN}
%                a\hei' = \SqCup_{x \in a\hei'} [x]
%            \end{EqLN}
%        \end{itemize}
%    }
%\end{lm}

%\subsection*{Central Derivations}
\begin{definition}
    $G$ is a stem group if $$ Z(G) \le G' $$
\end{definition}

\Bcr{hei:cor:3} can be generalised under the assumption that $G$ is a stem group. Note that $\hei$ is a stem group since $\hei' = Z(\hei)$ by \Bcr{heicommsub}. See \cite{Hall} for more details on (finite, which is not our case) stem groups.

\begin{pr}
    Let $G$ be a stem group. If $d$ is given by such $\chi$ that $\supp \chi \le \subgrd z$ for $z \in Z(G)$, then $d \in \Der_{G'}$.
\end{pr}

%Therefore
Central derivations were introduced in \cite{Aru20} as operators given on a group algebra generators $g\in G$ with the homomorphism $\tau: G\to (\cmplx ,+)$ and the central element $z\in Z(G)$ by formula
$$
    d_{\tau,z}: g\mapsto \tau(g) gz.
$$
Recall that \cite[Proposition~6]{Aru20} shows that central derivations form a Lie subalgebra in $\Der(G)$ (denote as $\ZDer(G)$).

\begin{pr}
     Let $G$ be NOT a stem group. Then there is an \textit{"induced"} nontrivial grading of $\ZDer(G)$ with $G'$.
\end{pr}

Note that if $G$ is a stem group, \textit{the grading is trivial}, i.e. not a grading at all in our terms. 

Another example of an induced grading follows.

    \innerdef
    
    Let $\IDer$ be the set of all derivations of the form ($a_1, \dots, a_n \in \cmplx; y_1, \dots, y_n \in G$)
    $$ G \ni x \mapsto a_1[x, y_1] + \dots + a_n[x,y_n] = [x, a_1y_1 + \dots + a_ny_n] $$

    A direct calculations shows that $\IDer$ is an ideal in $\Der$ (i.e. for any $d \in \IDer$ and for any $\de \in \Der:$ $[d, \de], [\de, d] \in \IDer$; recall that $[d, \de](x) = d(\de(x)) - \de(d(x))$.) 

    Therefore, there exists a factor-algebra $\ODer := \Der / \IDer$.
    
    \begin{corollary}
    If $G$ is NOT a stem group, there is an induced (non-trivial) grading of $\ODer$ with $G'$ of the form 
    
    $$ \ODer = \Oplus_{k \in G/G'} \Der_k / \IDer_k, $$

    where $\IDer_k := \Der_k \cap \IDer$
    
    \end{corollary}

%Induced gradings on central and outer and quasiouter derivations. The induced grading of central derivations is non-trivial if and only if G is not a stem group. G is a knot group (grading with Z)\todo{!!!}

%Get rid of Lemma4\todo{!!!}

%\begin{}

%However, there exists such $G$ such that the inverse is not true. Consider $G = GL_n$, then $Z(G) = ...$, $G' = SL_n$

%\bibliographystyle{plain} % We choose the "plain" reference style
%\bibliography{refs.bib} % Entries are in the refs.bib file

\begin{thebibliography}{9}
\bibitem{AruMisSht16}A. A. Arutyunov, A. Derivations of group algebras. {\em Fundam. Prikl. Mat.}. \textbf{21}, 65-78 (2016), http://mi.mathnet.ru/fpm1768
\bibitem{AleAru20}A. V. Alekseev, A. Derivations in semigroup algebras. {\em Eurasian Math. J.}. \textbf{11}, 9-18 (2020), http://mi.mathnet.ru/emj361
\bibitem{Aru20}Arutyunov, A. Derivation Algebra in Noncommutative Group Algebras. {\em Proceedings Of The Steklov Institute Of Mathematics}. \textbf{308}, 22-34 (2020,1), https://doi.org/10.1134/s0081543820010022
\bibitem{CREEDON2019247}Creedon, L. \& Hughes, K. Derivations on group algebras with coding theory applications. {\em Finite Fields And Their Applications}. \textbf{56} pp. 247-265 (2019), https://www.sciencedirect.com/science/article/pii/S107157971830145X
\bibitem{AruKos2021}Arutyunov, A. \& Kosolapov, L. Derivations of group rings for finite and FC groups. {\em Finite Fields And Their Applications}. \textbf{76} pp. 101921 (2021,12), https://doi.org/10.1016/j.ffa.2021.101921
\bibitem{kolesnikov2023prenovikov}Kolesnikov, P., Mashurov, F. \& Sartayev, B. On pre-Novikov algebras and derived Zinbiel variety.  (2023)
\bibitem{alekseev2020sigmatauderivations}Alekseev, A., Arutyunov, A. \& Silvestrov, S. On $(\sigma,\tau)$-derivations of group algebra as category characters.  (2020)
\bibitem{Hall}Hall, P. The classification of prime-power groups.. {\em Journal Für Die Reine Und Angewandte Mathematik}. \textbf{1940}, 130-141 (1940), https://doi.org/10.1515/crll.1940.182.130
\end{thebibliography}

\end{document}